\documentclass{article}
\usepackage{arxiv}
\usepackage{epsfig}
\usepackage{float}
\usepackage{amsmath,amssymb}
\usepackage{algorithm,algorithmic}
\usepackage[labelformat=simple]{subfig}

\usepackage{xcolor}
\usepackage{epstopdf}
\usepackage{multirow}
\usepackage{latexsym,bm}
\usepackage{doi}     

\graphicspath{{./fig/}}

\usepackage{hyperref}
\makeindex

\title{Reduced-order modeling for Ablowitz-Ladik equation}

\author{
  Murat Uzunca\\
   Department of Mathematics, Sinop University\\
     Sinop-Turkey\\ \
     \texttt{muzunca@sinop.edu.tr}\\
     \And
     B\"ulent Karas\"ozen\\
     Institute of Applied Mathematics \& Department of Mathematics\\
     Middle East Technical University\\
     Ankara-Turkey\\
     \texttt{bulent@metu.edu.tr}
}

\begin{document}

\maketitle

\begin{abstract}
In this paper, reduced-order models (ROMs) are constructed for the Ablowitz-Ladik equation (ALE), an integrable semi-discretization of the nonlinear Schr\"{o}dinger equation (NLSE) with and without damping. Both ALEs are non-canonical conservative and dissipative Hamiltonian systems with the Poisson matrix depending quadratically on the state variables, and with quadratic Hamiltonian. The full-order solutions are obtained with the energy preserving midpoint rule for the conservative ALE and exponential midpoint rule for the dissipative ALE. The reduced-order solutions are constructed intrusively by preserving the skew-symmetric structure of the reduced non-canonical Hamiltonian system by applying proper orthogonal decomposition (POD) with the Galerkin projection. For an efficient offline-online decomposition of the ROMs, the quadratic nonlinear terms of the Poisson matrix are approximated by the  discrete empirical interpolation method (DEIM). The computation of the reduced-order solutions is  further accelerated by the use of tensor techniques. Preservation of the Hamiltonian and momentum for the conservative ALE, and preservation of dissipation properties of the dissipative ALE, guarantee the long-term stability of soliton solutions.
\end{abstract}


\section{Introduction}

The nonlinear Schr\"{o}dinger equation (NLSE) is  a Hamiltonian partial differential equation (PDE)  modeling wave propagation phenomena in plasma physics, nonlinear optics, water waves, and in many other fields.
Ablowitz-Ladik equation (ALE) \cite{Ablowitz76} represents an integrable discretization of NLSE \cite{Islas01,Schober99,Tang07}.
It is not possible, in general, to preserve both the energy and the symplectic form of Hamiltonian systems at the same time.
Due to non-canonical Hamiltonian structure of the  ALE, the symplectic integrators are not directly applicable.
The energy preserving discrete gradient methods, such as 
the average vector field (AVF) method and the mid-point rule preserve the energy (Hamiltonian) and the quadratic invariants such as the momentum and norm of the ALE \cite{Cohen11}.
The linearly damped NLSE and the associated damped ALE are linearly perturbed non-canonical Hamiltonian systems. 
Discrete gradient methods  such as exponential mid-point rule have been also used for dissipative systems like the damped ALE \cite{Moore21}, which 
preserve the dissipation of energy, momentum, and norm with the correct rate.

Accurate and stable solution of PDEs with different discretization  techniques require large computing time and memory. In the last two decade,
many reduced-order models (ROMs) have been developed  as surrogate models.
By projecting the high fidelity full-order solutions of PDEs on low dimensional reduced spaces, the low-dimensional ROMs are constructed.
Among them, the proper orthogonal decomposition (POD) \cite{Berkooz93,Sirovich87}
is a widely used ROM technique. 
For an overview about model order reduction (MOR), we refer to the recent  book \cite{Benner20hbmor}.
Conventional ROMs do not ensure the conservation of the invariants like the energy. This may result in wrong or unstable reduced-order solutions.
ROMs that preserves energy-like structures are developed for Lagrangian systems \cite{Carlberg15} and port-Hamiltonian systems \cite{Chaturantabu16}. 
For canonical Hamiltonian systems like the linear wave equation, Sine-Gordon equation, NLSE, symplectic model reduction techniques ensuring long term stability of the reduced model, have been constructed in \cite{Peng16,Buchfink18,Hesthaven16,Buchfink20greedy}. 
The ROMs preserving the energy have been developed for non-canonical Hamiltonian PDEs like the Korteweg-de Vries (KdV) equation \cite{Gong17,Miyatake19,Karasozen21kdv}, NLSE \cite{Karasozen18nls} and rotating shallow water equations \cite{Karasozen21sw,Karasozen22rtswe}.  
All these structure-preserving ROMs are described for Hamiltonian systems in a recent review paper \cite{Hesthaven22an}.

Large number of reduced-order modes are required to conserve the invariants of the PDEs like the NLS, which contains cubic nonlinear terms. In the presence of nonlinearity, the online computational cost of the POD scales with the dimension of the full-order model. For a polynomial nonlinearity, the higher the degree of polynomial the more rapidly the cost of the tensorial POD increases, the cost scales with the reduced dimension to the power of the degree of the polynomial nonlinearity, which means that the online computational cost increases rapidly for the polynomial nonlinearity of degree three or higher.
For models depending only on quadratic nonlinearities, the tensorial POD, in general, has much lower computational cost than the POD \cite{Stefanescu14}. Therefore we consider the time transformed ALE equation, which contains quadratic and linear terms. We develop ROMs for the conservative and damped ALEs while preserving its physical properties like conservation/dissipation of the Hamiltonian and the momentum which is a quadratic invariant.
Following the approach in \cite{Gong17,Miyatake19} where the authors construct an energy preserving ROM for the KdV equation in non-canonical Hamiltonian form with constant skew-symmetric matrix, we construct ROMs that preserves the energy of the conservative ALE containing a state-dependent Poisson matrix.
Following the same approach with an additional linear damping term, a ROM is constructed that preserves the dissipation of the energy and the momentum of the dissipative ALE.
An important paradigm of MOR is the offline-online decomposition. 
The FOM and the POD basis are computed in the offline stage. In the online stage,  
the reduced system is solved. In order to separate the offline and online stage, the nonlinear term in the Poisson matrix is approximated with a hyper-reduction technique, i.e., discrete empirical interpolation method (DEIM) \cite{chaturantabut10nmr,gugercin16a}.  For dynamical systems such as the ALEs with wave-type solutions, a relatively large number of POD modes are needed to represent the physical behavior of the system in reduced-order form. Therefore, the reduced-order system should be solved efficiently. Here, the quadratic reduced system is solved utilizing tensor techniques \cite{Karasozen21sw,Benner21,Benner18} by the use of MULTIPROD \cite{leva08mmm} in order to speed up online computations.

Preservation of the energy conservation and dissipation by the ROMs for ALEs in non-canonical form is  challenging.  The main contributions  of the paper are:

\begin{itemize}

\item 
The reduced Hamiltonian and momentum of the conservative ALE is preserved, leading to accurate and stable soliton solutions.

\item Dissipation of the reduced Hamiltonian and momentum of the linearly perturbed ALE is preserved with the correct rate by the ROMs.

\item The POD/DEIM reduced systems are solved through a tensorial framework, which enables computationally efficient offline-online decomposition.

\end{itemize}

The paper is organized as follows. In Section \ref{fom}, the high fidelity discretization that preserves the conservation/dissipation property of the conservative/damped ALEs is presented. In Section \ref{rom}, ROMs for both types of ALEs are described. Two numerical tests in Section \ref{num} verify the  conservation and dissipation of the energy/momentum by the  ROMs. Conclusions are given  in Section \ref{con}.

\section{Discretization of the conservative and dissipative Ablowitz-Ladik equations}
\label{fom}

The NLSE  is one of the most known  nonlinear PDE in nonlinear optics, quantum physics, plasma physics.
The NLSE is given as
\begin{equation} \label{nls}
i \dot{\psi}  =  -\psi_{xx} - 2\gamma |\psi|^2\psi,
\end{equation}
where $\psi(x,t)$ denotes the complex-valued wave function. We impose
 periodic boundary conditions, $\psi(-L,t) = \psi(L,t)$, for $t \in (0,T]$ with a final time $T>0$. For $\gamma > 0$, the NLSE \eqref{nls} possesses infinitely many integrals such as energy, momentum and norm \cite{Tang07,Zaharov74}, i.e., it is completely integrable.

ALE \cite{Ablowitz76} represents an integrable Hamiltonian semi-discretization of the NLSE \eqref{nls} \cite{Islas01,Schober99,Tang07}
\begin{equation} \label{al1}
i\dot{\psi}_n = - \frac{1}{h^{2}}\left(\psi_{n+1}-2\psi_{n}+\psi_{n-1}\right) - \gamma |\psi_n|^2\left(\psi_{n+1}+\psi_{n-1}\right) =0,
\end{equation}
with
$h=2L/N$, $\psi_{n}=\psi (x_n,t)$, $x_n=-L+(n-1)h$, $n= 1, \ldots, N+1 $. The solutions of the ALE \eqref{al1} converge to the solutions of the NLSE \eqref{nls} when the step-size $h\rightarrow 0$. Under unitary time dependent transformation
$\psi_n \rightarrow w_n e^{-2i t/h^2}$, the ALE \eqref{al1} becomes  a non-canonical Hamiltonian system
\cite{Islas01,Tang07}
\begin{equation}\label{al2}
i \dot{w}_n = - \frac{1}{h^2}(w_{n+1} + w_{n-1}) \left ( 1  + \gamma h^2|w_n|^2 \right ).
\end{equation}
Separating the real and complex parts as
$w = p + i q$, the equation \eqref{al2} yields the coupled system
\begin{equation}\label{als}
\begin{aligned}
\dot{p}_n&=-\frac{1}{h^{2}}\left(1+\gamma h^2\left(p_{n}^2+ q_n^2\right)\right)\left(q_{n+1}+q_{n-1}\right),  \\
\dot{q}_n&=\frac{1}{h^{2}}\left(1+\gamma h^2\left(p_{n}^2+ q_n^2\right)\right)\left(p_{n+1}+p_{n-1}\right).
\end{aligned}
\end{equation}
Setting $\bm{p}=(p_1,\ldots , p_N)^T$ and $\bm{q}=(q_1,\ldots , q_N)^T$, the  ALE \eqref{als} is given in matrix-vector form by
\begin{align}\label{alnc}
\begin{pmatrix}
\dot{\bm{p}}\\\dot{\bm{q}}
\end{pmatrix}=
\begin{pmatrix}
0&-M(\bm{p},\bm{q})\\
M(\bm{p},\bm{q})&0
\end{pmatrix}
\begin{pmatrix}
\nabla_p H(\bm{p},\bm{q})\\
\nabla_  q H(\bm{p},\bm{q})
\end{pmatrix},
\end{align}
with the quadratic Hamiltonian  $H(\bm{p},\bm{q})$,
 and the matrix $M(\bm{p},\bm{q})$ with nonlinear terms given by
\begin{equation} \label{ham1}
\begin{aligned}
 H(\bm{p},\bm{q})&=\frac{1}{h^{2}}\sum_{n=1}^{N}\left(p_{n}p_{n-1}+q_{n}q_{n-1}\right), \\
 M(\bm{p},\bm{q})&= \text{diag}\left(d_1,\ldots ,d_N\right) , \qquad d_n = 1+\gamma h^{2} \left(p^2_n+q^2_n\right).
\end{aligned}
\end{equation}
Under the above setting, the gradients of the Hamiltonian reduce to
\begin{equation} \label{nabla}
\nabla_p H(\bm{p},\bm{q})  = D {\bm p}, \quad \nabla_q H(\bm{p},\bm{q})  = D {\bm q},\quad
D =\frac{1}{h^{2}}
\begin{pmatrix}
0 & 1 & &  & 1\\
1 & 0 & 1& & \\
 & \ddots & \ddots & \ddots & \\
1 & &  & 1 & 0
\end{pmatrix}.
\end{equation}

Besides the  Hamiltonian, there exist quadratic invariants, i.e., Casimirs, such as the momentum ${\mathcal I}$ whose discrete form is given by
\begin{equation}\label{mom1}
 I(\bm{p},\bm{q}) = \sum_{n=1}^N \left ( q_n\frac{p_{n+1} -p_{n-1}}{2h} - p_n\frac{q_{n+1} -q_{n-1}}{2h}     \right ).
\end{equation}

For the solution vector $\bm z= (\bm p^T, \bm q^T)^T$, the ALE \eqref{alnc} is a $2N$-dimensional skew-gradient ODE of the form
\begin{equation} \label{sgfom}
	\dot{\bm{z}} = S(\bm{z})\nabla_zH(\bm{z}),
\end{equation}
where the  skew-symmetric matrix $S(\bm{z})$ is given by
$$
S(\bm{z})=\begin{pmatrix}
0&-M(\bm{p},\bm{q})\\
M(\bm{p},\bm{q})&0
\end{pmatrix}.
$$
The AVF method applied to the skew-gradient ODE system \eqref{sgfom} reads as
\begin{equation} \label{alncavf}
	\frac{\bm{z}^{k+1}-\bm{z}^{k}}{\Delta t} = S\left( \frac{\bm{z}^{k+1}+\bm{z}^{k}}{2} \right)\int_0^1 \nabla_zH(\xi\bm{z}^{k+1}+(1-\xi)\bm{z}^{k})d\xi, \quad k=1,\ldots,K,
\end{equation}
where $\bm{z}^{k}=\bm{z}(t_k)$.
For quadratic Hamiltonians as for the ALE \eqref{alnc}, the AVF method is equivalent to the midpoint rule	
\begin{equation} \label{almp}
	\frac{\bm{z}^{k+1}-\bm{z}^{k}}{\Delta t} = S\left( \frac{\bm{z}^{k+1}+\bm{z}^{k}}{2} \right)\nabla_zH\left( \frac{\bm{z}^{k+1}+\bm{z}^{k}}{2} \right),
\quad k=1,\ldots,K.
\end{equation}

The damped NLSE
\begin{equation} \label{dnls}
i \dot{\psi}  =  -\psi_{xx} - 2\gamma |\psi|^2\psi - i \mu\psi,
\end{equation}
with the damping factor $\mu$, describes resonant
phenomena in nonlinear media  \cite{Fu16}. 
The Ablowitz-Ladik discretization of the damped NLSE \eqref{dnls} yields
\begin{equation} \label{aldamped}
\begin{aligned}
\dot{p}_n&=-\frac{1}{h^{2}}\left(1+\gamma h^2\left(p_{n}^2+ q_n^2\right)\right)\left(q_{n+1}+q_{n-1}\right) -\mu p_n,  \\
\dot{q}_n&=\frac{1}{h^{2}}\left(1+\gamma h^2\left(p_{n}^2+ q_n^2\right)\right)\left(p_{n+1}+p_{n-1}\right) - \mu q_n,
\end{aligned}
\end{equation}
which represent a linearly perturbed non-canonical Hamiltonian system \cite{Moore21}.
In matrix-vector form, the system \eqref{aldamped} is a skew-gradient system with a damping term
\begin{equation} \label{sgfomdamped}
	\dot{\bm{z}} = S(\bm{z})\nabla_zH(\bm{z}) - \mu{\bm z}.
\end{equation}
The energy (Hamiltonian) and the momentum  dissipate like
$$
\frac{d}{dt} {\mathcal I} = -\mu {\mathcal I}\; \Rightarrow \; {\mathcal I} (\bm z (t) ) = e^{-\mu t} {\mathcal I} ({\bm z}^0).
$$
Dissipation of the energy is expressed in form of the energy balance equation \cite{Moore21} given by
\begin{equation} \label{hambal}
H(e^{X_1}\bm{z}^{k+1}) - H(e^{X_0}\bm{z}^{k}) = 0, \quad R_H = \ln \left( \frac{H(\bm{z}^{k+1})}{H(\bm{z}^{k})}\right) + \Delta t \mu,
\end{equation}
where
$$
X_{\alpha} = (t_{k+\alpha} - t_{k+ \frac{1}{2}})\mu, \quad X_0 = -\frac{\Delta t}{2}\mu, \quad X_1 =  \frac{\Delta t}{2}\mu.
$$
Similarly,  the dissipation of the momentum  \eqref{mom1} is given as \cite{Moore21}
\begin{equation} \label{mombal}
e^{2X_1}I(\bm{z}^{k+1}) - e^{2X_0}I(\bm{z}^{k}) = 0, \quad
R_I = \ln \left( \frac{ I (\bm{z}^{k+1})}{ I (\bm{z}^{k})}\right) - 2\Delta t \mu.
\end{equation}

The exponential mid-point rule preserves the correct dissipation rate of the energy and dissipative dynamics of the quadratic  Casimirs such as the momentum and norm \cite{Moore21}.
Approximation of   the gradient of the Hamiltonian with
\begin{equation}
\bar{\nabla} H (e^{X_1} {\bm z}^{k+1} , e^{X_0}  {\bm z}^k )
=\int_0^1 \nabla H (\eta e^{X_1} {\bm z}^{k+1} + (1-\eta) e^{X_0}{\bm z}^k d\eta,
\end{equation}
yields the exponential midpoint rule  \cite{Moore21,Bhatt16}
\begin{equation} \label{exmd}
\frac{e^{X_1} {\bm z}^{k+1} - e^{X_0}  {\bm z}^k  }{\Delta t} =
S\left (   \frac{e^{X_1} {\bm z}^{k+1} + e^{X_0}  {\bm z}^k  }{2} \right ) \bar{\nabla} H (e^{X_1} {\bm z}^{k+1} , e^{X_0}  {\bm z}^k ).
\quad k=1,\ldots,K,
\end{equation}

The exponential mid-point rule \eqref{exmd} preserves the dissipation rate of the invariants of the damped ALE \eqref{aldamped}
$$
{\mathcal I}^{k+1} = e^{-\mu\Delta t  } {\mathcal I}^{k}, \; k=1,\ldots,K.
$$

Both the midpoint rule \eqref{almp} and the exponential midpoint rule \eqref{exmd} 
are time-reversible, therefore second order convergent and unconditionally stable time integrators.

\section{Reduced-order modelling}
\label{rom}

The ROM for the damped ALE differs only by adding the linear damping term, therefore we describe the construction of the ROMs for the conservative ALE \eqref{alnc} whose compact form is given by \eqref{sgfom}.

\subsection{POD reduced skew-gradient system}

The ROMs are constructed by projecting the full-order system onto a low dimensional reduced space spanned by POD basis.
Let ${\mathsf S}_p$ and ${\mathsf S}_q$ denote the matrix of solution snapshots given by
$$
{\mathcal S}_p = \left[{\bm p}^1   \cdots  {\bm p}^{K} \right] \in\mathbb{R}^{N\times K}, \qquad
{\mathcal  S}_q = \left[{\bm q}^1   \cdots  {\bm q}^{K} \right]   \in\mathbb{R}^{N\times K},
$$
where ${\bm p}^k, {\bm q}^k\in\mathbb{R}^N$, $k=1,\ldots ,K$, are the fully discrete solution vectors of the full-order ALE \eqref{alnc} through the mid-point rule \eqref{almp}.
Then, the POD basis modes are taken as the left singular vectors related to the most dominant singular values
of the snapshot matrices ${\mathsf S}_p$ and  ${\mathsf S}_q$
\begin{equation*}
{\mathcal S }_p= W_p \Sigma_p Z_p^T\, \quad {\mathcal S }_q= W_q \Sigma_q Z_q^T,
\end{equation*}
where for $*\in\{p,q\}$, $ W_{*} \in \mathbb{R}^{ N\times K}$ and
$Z_{*}\in \mathbb{R}^{ K\times K}$ are orthonormal matrices whose column vectors are the left and right singular vectors, respectively, and $\Sigma_{*}\in \mathbb{R}^{ K\times K}$
is the diagonal matrix containing the singular values $\sigma_{*,1} \ge \sigma_{*,2} \ge \cdots \ge \sigma_{*,K}\geq 0$.
For some positive integer $N_r\ll \min\{N,K\}$, let $V_{*}\in\mathbb{R}^{N\times N_r}$ denotes the truncated matrix of POD modes consisting of the first $N_r$ left singular vectors from $W_{*}$. In this paper, we take the same number of POD modes for $p$ and $q$. 
The POD matrix $V_{*}$  is the minimizer of the least squares error
$$
\min_{V_{*}\in \mathbb{R}^{ N\times N_r}}||{\mathsf S}_{*}-V_{*}V_{*}^T{\mathsf S}_{*} ||_F^2 = \sum_{j=N_r+1}^{K} \sigma_{*,j}^2,
$$
with $\|\cdot\|_F$ denoting the Frobenius norm.
Once the POD basis is obtained, an approximation for the full-order solutions of \eqref{alnc} from the reduced space spanned by the POD modes, is given  as
\begin{equation}\label{relz}
{\bm p} \approx  \widehat{{\bm p}}=\ V_{p}{\bm p}_r, \quad {\bm q} \approx  \widehat{{\bm q}}= V_{q}{\bm q}_r,
\end{equation}
where ${\bm p}_r,{\bm q}_r : [0,T]\mapsto \mathbb{R}^{N_r}$ are the coefficient vectors. The coefficient vectors are the solutions of the $2N_r$-dimensional reduced system
\begin{align} \label{galpod1}
\dot{{\bm z}}_r = V_{z}^T S(\widehat{{\bm z}})\nabla_{\bm z} H(\widehat{{\bm z}}),
\end{align}
which is constructed by the Galerkin projection onto the reduced space. In the ROM \eqref{galpod1},
${\bm z}_r=({\bm p}_r^T,{\bm q}_r^T)^T: [0,T]\mapsto \mathbb{R}^{2N_r}$ is the vector of coefficients,
$\widehat{{\bm z}}= V_{z}{\bm z}_r=(\widehat{{\bm p}}^T,\widehat{{\bm q}}^T)^T : [0,T]\mapsto \mathbb{R}^{2N} $ is the vector of reduced approximations, and
the diagonal matrix $V_{z}$ contains the POD matrices to the state variables
\begin{equation*}
V_{z}=
\begin{pmatrix}
V_{p} & \\
& V_{q}
\end{pmatrix}\in\mathbb{R}^{2N\times 2N_r}.
\end{equation*}

Unlike the FOM \eqref{sgfom}, the reduced system \eqref{galpod1} is not in skew-gradient form. To recover a reduced system in skew-gradient form,
we  formally insert $V_{z}V_{z}^T\in \mathbb{R}^{2N \times 2N}$ between $S(\widehat{{\bm z}})$ and $\nabla_{\bm z} H(\widehat{{\bm z}})$ \cite{Gong17,Miyatake19}. This results in the POD reduced skew-gradient system
\begin{align} \label{sgrom}
\dot{{\bm z}}_r = S_r({\bm{z}_r})\nabla_{\bm{z}_r} \widetilde{H}({\bm{z}_r}),
\end{align}
where $S_r({\bm{z}_r})= V_{z}^T S( V_z\bm{z}_r)V_{z}\in\mathbb{R}^{2N_r\times 2N_r}$ is the reduced skew-symmetric matrix, $\nabla_{\bm{z}_r} \widetilde{H}({\bm{z}_r})= V_{z}^T\nabla_{\bm{z}}H(V_{z}{\bm{z}_r})\in\mathbb{R}^{2N_r}$ is the reduced discrete gradient of the Hamiltonian, and $\widetilde{H}({\bm{z}_r})=H(V_z{\bm{z}_r})$.
Then, the skew-gradient structure of the
reduced system \eqref{sgrom} yields
\begin{align*}
\frac{d}{dt}H(\widehat{{\bm z}}) = \frac{d}{dt}H(V_z{\bm z}_r)&= \left[V_z^T\nabla_{\mathbf{z}} H(V_z{\bm z}_r)\right]^T\dot{{\bm z}}_r\\
&= \left[\nabla_{\mathbf{z}_r} \widetilde{H}({\bm z}_r)\right]^TS_r({\bm z}_r)\left[\nabla_{\mathbf{z}_r} \widetilde{H}({\bm z}_r)\right]= 0,
\end{align*}
which means that the reduced Hamiltonian $\widetilde{H}({\bm z}_r)$ is preserved by  \eqref{sgrom}.

From the preservation of the skew-symmetry of the reduced system \eqref{sgrom}, the energy dissipation is also preserved for dissipative systems such as the damped ALE \cite{Gong17}.
Both the reduced systems are solved by the same time integrators as for the FOMs, i.e., the reduced ALE \eqref{sgrom} is solved in time with the implicit midpoint rule \eqref{almp}, whereas the reduced dissipative ALE \eqref{sgfomdamped} is solved with the exponential midpoint rule \eqref{exmd}.

\subsection{POD-DEIM ROM}

Explicit form of the POD reduced skew-gradient system \eqref{sgrom} reads as
\begin{equation} \label{alred1}
\begin{pmatrix}
\dot{{\bm p}}_r\\ \dot{{\bm q}}_r
\end{pmatrix} =
\begin{pmatrix}
0&- V_p^T M(V_z\bm{z}_r)V_q \\
V_q^T M(V_z\bm{z}_r)V_p & 0
\end{pmatrix}
\begin{pmatrix}
V_p^T\nabla_p H (V_z\bm{z}_r)  \\
V_q^T\nabla_q H (V_z\bm{z}_r)
\end{pmatrix},
\end{equation}
where $M(V_z\bm{z}_r)=M(V_p\bm{p}_r,V_q\bm{q}_r)=M(\widehat{\bm{p}},\widehat{\bm{q}})$ is given as in \eqref{ham1}, i.e., $M(V_z\bm{z}_r)=\text{diag}\left(d_1,\ldots ,d_N\right)$ with $d_n = 1 + \gamma h^{2} \left(\widehat{p}^2_i+\widehat{q}^2_i\right)$, $n=1,\ldots,N$.
Using the diagonal structure of the matrix $M(V_z\bm{z}_r)$, setting the nonlinear vector ${\bm m} =  (d_1,\ldots ,d_N)^T$, and inserting the discrete gradients of the Hamiltonian in \eqref{nabla}, i.e., $\nabla_p H (\bm{z})=D\bm{p}$ and $\nabla_q H (\bm{z})=D\bm{q}$, \eqref{alred1} can be rewritten as 
\begin{equation}\label{alred2}
\begin{aligned}
\dot{{\bm p}}_r &= -V_p^T\left[ \bm{m}\odot (V_qV_q^TDV_q{\bm q}_r)  \right], \\
\dot{{\bm q}}_r &= V_q^T\left[ \bm{m}\odot (V_pV_p^TDV_p{\bm p}_r)  \right],
\end{aligned}
\end{equation}
with $\odot$ denoting the element-wise product of vectors.
Due to the nonlinear function ${\bm m} =  (d_1,\ldots ,d_N)^T: [0,T] \mapsto \mathbb{R}^{N}$, the reduced system  \eqref{alred2} still depends on the FOM dimension. This can be circumvent applying the DEIM \cite{chaturantabut10nmr,gugercin16a}.  The idea in DEIM is to interpolate the nonlinear vector ${\bm m}$ using only $N_d\ll\min\{N,K\}$ entries. For this goal,
one needs to compute the interpolation basis and an operator which selects the interpolation points and calculate the empirical basis. Let $\Phi = [\phi_{1}, \dots, \phi_{N_d} ]\in \mathbb{R}^{N \times N_d}$ denotes the $N_d$-dimensional interpolation basis.  The interpolation basis $\{\phi_1,\dots,\phi_{N_d}\}$ is constructed by applying the POD to the snapshot matrix $\mathsf{S}_{m}$ containing the  nonlinear terms, given by
\begin{equation}
\mathsf{S}_{m} = [ {\bm m}^1 \cdots  {\bm m}^{K}]\in\mathbb{R}^{N\times K},
\label{nonsvd}
\end{equation}
where ${\bm m}^k={\bm m }(t_k)$ denotes the nonlinear vector computed by the fully discrete solutions of the FOM \eqref{alnc} at $t=t_k$.
In other words, the DEIM basis modes $\phi_i\in \mathbb{R}^{N }$ are determined as the ${N_d}$ left singular vectors related to the first ${N_d}$ largest singular values
of the snapshot matrix $\mathsf{S}_m$. Once the DEIM basis $\Phi$ is constructed, the nonlinearity can be approximated as
\begin{equation}\label{deimapp}
{\bm m}(t) \approx \Phi \mathsf{c}(t),
\end{equation}
where $\mathsf{c}(t) : [0,T] \mapsto \mathbb{R}^{N_d}$. 
In order to uniquely solve for the coefficient vector $\mathsf{c}$, the overdetermined system \eqref{deimapp} needs to be projected by multiplication from left by a matrix, say $\mathsf{P}$. 
The matrix $\mathsf{P}$ is indeed a selection matrix which is determined by a greedy algorithm \cite{chaturantabut10nmr}.
Alternatively, the Q-DEIM \cite{gugercin16a} is also used to determine the matrix $\mathsf{P}$, which uses the pivoted QR-factorization, and it is more accurate and stable than the DEIM, see Algorithm~\ref{alg:qdeim}.

\begin{algorithm}[h]
\caption{\label{alg:qdeim}Q-DEIM algorithm.}
\begin{algorithmic}[1]
	\STATE\textbf{Input:} Basis matrix $\Phi\in\mathbb{R}^{N\times p}$
	\STATE\textbf{Output:} Selection matrix $\mathsf{P}$
	\STATE Perform pivoted QR factorization of $\Phi^{\top}$ so that $\Phi^{\top}\Pi = QR$
	\STATE Set $\mathsf{P} = \Pi (:,1:p)$
\end{algorithmic}
\end{algorithm}

After computation of the selection matrix $\mathsf{P}$,
the coefficient vector $\mathsf{c}(t)$ is computed by solving the projected linear system
$\mathsf{P}^{\top}\Phi \mathsf{c}(t) = \mathsf{P}^{\top}{\bm m}(t)$,
by which the approximation \eqref{deimapp} becomes
\begin{equation} \label{fkapprox}
{\bm m}(t) \approx \Psi{\bm m}_r (t),
\end{equation}
where 
 $\Psi:=\Phi(\mathsf{P}^{\top}\Phi)^{-1}\in \mathbb{R}^{ N\times {N_d}}$ is a constant matrix which is precomputed in the offline stage. The reduced nonlinear vector ${\bm m}_r (t):=\mathsf{P}^{\top} {\bm m} (t):[0,T]\mapsto \mathbb{R}^{N_d} $, is computed in the online stage with the reduced dimension ${N_d}\ll N$. In fact, the reduced nonlinear vector ${\bm m}_r (t)$ is nothing but $N_d$ selected entries among the $N$ entries of the nonlinear vector ${\bm m}$.

Inserting the DEIM approximation ${\bm m} \approx \Psi{\bm m}_r$ into \eqref{alred2}, we obtain the skew-gradient POD-DEIM ROM	
\begin{equation}\label{deim}
\begin{aligned}
\dot{{\bm p}}_r &= -V_p^T\left[ (\Psi{\bm m}_r)\odot (V_qV_q^TDV_q{\bm q}_r)  \right], \\
\dot{{\bm q}}_r &= V_q^T\left[ (\Psi{\bm m}_r)\odot (V_pV_p^TDV_p{\bm p}_r)  \right].
\end{aligned}
\end{equation}

Because the nonlinearities in the ROM are approximated by the DEIM,
the Hamiltonian is preserved approximately. However, an upper bound for the preservation of the discrete energy by POD-DEIM reduced system can be derived similar to the derivations in \cite{Karasozen22rtswe,Karasozen18}, by the use of well-known DEIM approximation bound \cite[Lemma 3.2]{chaturantabut10nmr}
$$
\|{\bm m}-\Psi{\bm m}_r\|_2 \le \|(\mathsf{P}^{\top}\Phi)^{-1}\|_2\|(I-\Phi\Phi^{\top}){\bm m}\|_2,
$$
where $\|\cdot\|_2$ denotes the Euclidean $2$-norm, and $I$ is the identity matrix.

\subsection{The reduced quadratic system}

For the POD-DEIM ROM \eqref{deim}, the online and offline stages are still not separated.
The precomputable constant matrices and the reduced terms can be separated by the use of tensor techniques, so that the solution  of the POD-DEIM reduced skew-gradient system \eqref{deim} is accelerated.  We first rewrite the reduced system \eqref{deim} in which Kronecker product $\otimes$ takes place of component-wise product $\odot$. This can be handled by using a so-called matricized tensor $G\in\mathbb{R}^{N\times N^2}$ satisfying $G(\bm{a}\otimes\bm{b})=\bm{a}\odot\bm{b}$, i.e., by the left multiplication by the matrix $G$, we select the appropriate entries of the Kronecker product to get component-wise product. In terms of Kronecker product, the reduced system \eqref{deim} yields
\begin{equation}\label{deim2}
\begin{aligned}
\dot{{\bm p}}_r &= -V_p^TG\left[ (\Psi{\bm m}_r)\otimes (V_qV_q^TDV_q{\bm q}_r)  \right], \\
\dot{{\bm q}}_r &= V_q^TG\left[ (\Psi{\bm m}_r)\otimes (V_pV_p^TDV_p{\bm p}_r)  \right].
\end{aligned}
\end{equation}
Using the properties of Kronecker product, the system \eqref{deim2} further reduces to
\begin{equation}\label{deim3}
\begin{aligned}
\dot{{\bm p}}_r &= -V_p^TG\left( \Psi\otimes (V_qV_q^TDV_q) \right)\left( {\bm m}_r\otimes {\bm q}_r  \right), \\
\dot{{\bm q}}_r &= V_q^TG\left( \Psi\otimes (V_pV_p^TDV_p)  \right)\left( {\bm m}_r\otimes {\bm p}_r  \right),
\end{aligned}
\end{equation}
where it needs only the computation of the nonlinear vectors ${\bm m}_r\otimes {\bm p}_r:[0,T]\mapsto\mathbb{R}^{N_rN_d}$ and ${\bm m}_r\otimes {\bm q}_r:[0,T]\mapsto\mathbb{R}^{N_rN_d}$ in the online stage, and the terms  $V_p^TG\left( \Psi\otimes (V_qV_q^TDV_q) \right)\in\mathbb{R}^{N_r\times (N_rN_d)}$ and $V_q^TG\left( \Psi\otimes (V_pV_p^TDV_p)  \right)\in\mathbb{R}^{N_r\times (N_rN_d)}$ are constant matrices to be computed in the offline stage.
By DEIM approximation with tensor setting, online computations scale with $\mathcal{O}(N_dN_r^2)$, whereas it scales with $\mathcal{O}(N_rN^2)$ if DEIM approximation is not used.

Apart from the computational efficiency in the online stage, we also follow a computationally efficient approach in the offline stage for the calculation of the constant matrices $G\left( \Psi\otimes (V_qV_q^TDV_q) \right)\in\mathbb{R}^{N\times (N_rN_d)}$ and $G\left( \Psi\otimes (V_pV_p^TDV_p)  \right)\in\mathbb{R}^{N\times (N_rN_d)}$.
Since  the matrix $\widehat{G}:=G\left( \Psi\otimes (V_qV_q^TDV_q) \right)$  scales with the dimension of the FOM, the explicit calculation of  $\widehat{G}$ is inefficient. Using the structure of $\widehat{G}$, the matrix $\widehat{G}$ is given in MATLAB notation without constructing the matricized tensor $G$ explicitly by
\begin{align}\label{goyal}
	\widehat{G} =
	\begin{pmatrix}
	\Psi (1,:)\otimes (V_qV_q^TDV_q)(1,:)\\
	\vdots\\
	\Psi(N,:)\otimes (V_qV_q^TDV_q)(N,:)
	\end{pmatrix}.
	\end{align}
However, the computation in \eqref{goyal} needs $N$ for loops in which a matrix product is done. This drawback can be overcome in a computationally efficient way by the use of MULTIPROD \cite{leva08mmm}.
More clearly, by the properties of Kronecker product, the $i$th row of the matrix $\widehat{G}$ in \eqref{goyal} is given equivalently by
\begin{equation}\label{goyal2}
\widehat{G}(i,:)=(\text{vec}(\Psi(i,:)^\top (V_qV_q^TDV_q)(i,:))^\top , \quad i=1,\ldots,N.
\end{equation}
For each $i=1,\ldots,N$, all the operations in \eqref{goyal2} can be done at once by MULTIPROD as follows: we reshape the two-dimensional array (matrix) $\Psi\in \mathbb{R}^{N\times N_d}$ as a three-dimensional array $\widetilde{\Psi} \in \mathbb{R}^{N\times 1 \times N_d}$, and then we compute MULTIPROD of $(V_qV_q^TDV_q)\in \mathbb{R}^{N\times N_r \times 1}$ and $\widetilde{\Psi}\in \mathbb{R}^{N\times 1 \times N_d}$ in $2$nd and $3$th dimensions, which results in the three-dimensional array
\begin{equation}\label{goyal3}
\mathcal{\widehat{G}} =\text{MULTIPROD}( V_qV_q^TDV_q,\widetilde{\Psi})\in \mathbb{R}^{N\times N_r \times N_d},
\end{equation}
where the tensor $\mathcal{\widehat{G}}$ collects all the matrix product of two matrices of sizes $N_r\times 1$ and $1\times N_d$ within $N$ iterations.
Finally, the required matrix $\widehat{G}\in \mathbb{R}^{N\times (N_r N_d)}$ is obtained by reshaping the tensor $\mathcal{\widehat{G}}$ in \eqref{goyal3} into a two-dimensional array of dimension $N\times (N_r N_d)$.
Hence, utilizing the MULTIPROD, all the matrix products are calculated  in a single loop, that yields fast computations in the offline stage.

\section{Numerical results}
\label{num}

In this section, we illustrate the energy and momentum preserving properties for the conservative ALE and preservation of the dissipative structure of the damped ALE.
The relative cumulative energy criterion
\begin{equation} \label{energy_criteria}
\min_{1\leq p \leq K}\frac{\sum_{j=1}^p \sigma_{j}^2}{\sum_{j=1}^{K} \sigma_{j}^2  } > 1 - \kappa,
\end{equation}
is used to determine the number of POD/DEIM modes with a tolerance $\kappa$ , and with $p=N_r$ or $p=N_d$.

To asses the accuracy of the ROMs, we use the relative $L_2$-norm errors  between full- and reduced-order solutions
\begin{align}\label{relerr}
\|\bm{\psi}-\widehat{\bm  \psi}\|_{rel}=\frac{1}{K}\sum_{k=1}^{K}\frac{\|{\bm  \psi}^k-\widehat{\bm  \psi}^k\|_{L^2}}{\|{\bm  \psi}^k\|_{L^2}}.
\end{align}

Conservation of the discrete energy \eqref{ham1} and discrete momentum \eqref{mom1} of the FOMs and ROMs for the conservative ALE \eqref{al1} are  measured using the relative errors 
\begin{equation} \label{conserr1}
\|H\|_{\text{abs}} = \frac{1}{K}\sum_{k=1}^{K} \frac{ |H(\bm{w}^k)-H(\bm{w}^0)|}{ |H(\bm{w}^0)|  }, \quad
\|I\|_{\text{abs}} = \frac{1}{K}\sum_{k=1}^{K} \frac{ |I(\bm{w}^k)-I(\bm{w}^0)|}{ |I(\bm{w}^0)|  }.
\end{equation}

Similarly the dissipation rate of the discrete energy \eqref{hambal} and the discrete momentum \eqref{mombal} of the damped ALE \eqref{aldamped} are measured using the relative errors 
\begin{equation} \label{conserr2}
\|R_H\|_{\text{abs}} = \frac{1}{K}\sum_{k=1}^{K} \frac{ |R_H(\bm{w}^k)-R_H(\bm{w}^0)|}{ |R_H(\bm{w}^0)|  }, \quad
\|R_I\|_{\text{abs}} = \frac{1}{K}\sum_{k=1}^{K} \frac{ |R_I(\bm{w}^k)-R_I(\bm{w}^0)|}{ |R_I(\bm{w}^0)|  }.
\end{equation}

In all the simulations, unless otherwise stated, the ROM refers to the POD-DEIM ROM with tensorial setting. At the end, we present CPU time plots of different ROMs, i.e., POD ROM and POD-DEIM ROM with tensorial setting, in order to compare the computational efficiency between the ROMs with and without DEIM approximation.

\subsection{Conservative ALE}

\begin{figure}[H]
	\centerline{\includegraphics[width=0.5\columnwidth]{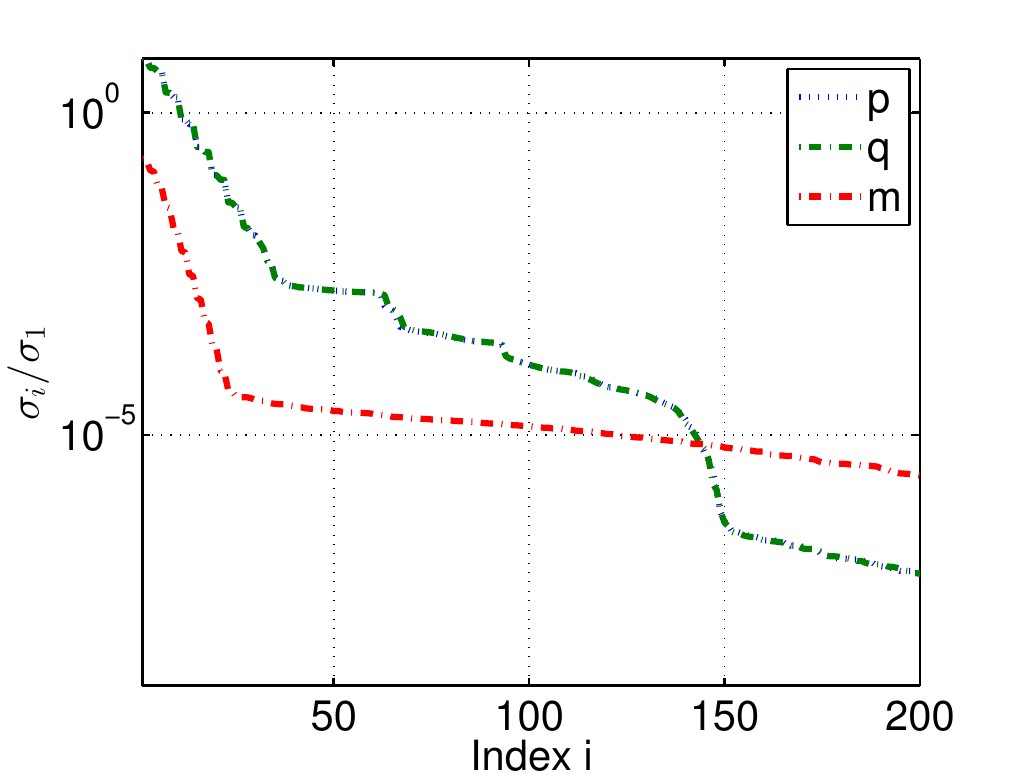}}
\caption{Singular values to the states and the nonlinearity}
\label{singal}
\end{figure}

We consider the ALE \eqref{al1} with soliton solution for $\gamma =1$ with the  initial condition \cite{Schober99}
$$
\psi(x,0) = 2\eta e^{2i\xi x_n} {\mathrm sech}(2\eta x_n),
$$
where $x_n=-50+0.5(n-1), n=1,\ldots N$, with $N=200$, $h=0.5$, $\eta =0.05$, $\xi=0.5$, and $\Delta t =0.01$ for $0 < t\le 50$.
The snapshot matrices are computed by saving the full discrete solutions
$\bm p$ and $\bm q$ at every five time steps, and they are of size $200\times 1000$. Figure~\ref{singal} shows the decay of singular values of related snapshot matrices. 
The singular values of the snapshots corresponding to states and nonlinear terms decay rapidly at the beginning, then they decay slowly in Figure \ref{singal}, 
which is characteristic for wave type and transport dominated problems  \cite{Ahmed20,Ohlberger16}.
The rate of the singular value decay is related to the Kolmogorov $r-$width describing the
error of a projection onto the best-possible space of a given dimension $r$.
It determines the linear reducibility of the FOM \cite{Peherstorfer22}. 
Therefore, relatively large number of POD/DEIM modes are needed to preserve the reduced conserved quantities. In the numerical experiments, we set the tolerances $\kappa = 10^{-4}$ and $\kappa = 10^{-6}$ in \eqref{energy_criteria}, giving  respectively $N_r=21$ POD and $N_d=21$ DEIM modes. 
The relative FOM-ROM error \eqref{relerr} is 9.55e-03, whereas, the relative energy and momentum errors are  1.90e-04 and 2.29e-04, respectively.

\begin{figure}[H]
	\centerline{\includegraphics[width=0.5\columnwidth]{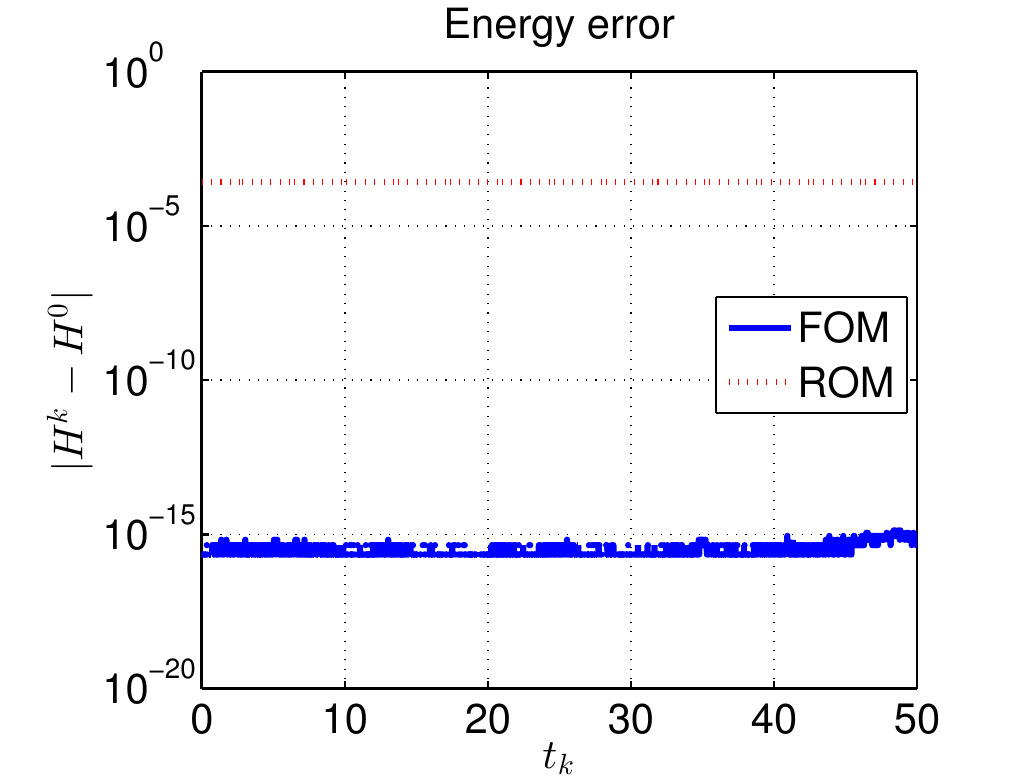}
		 \includegraphics[width=0.5\columnwidth]{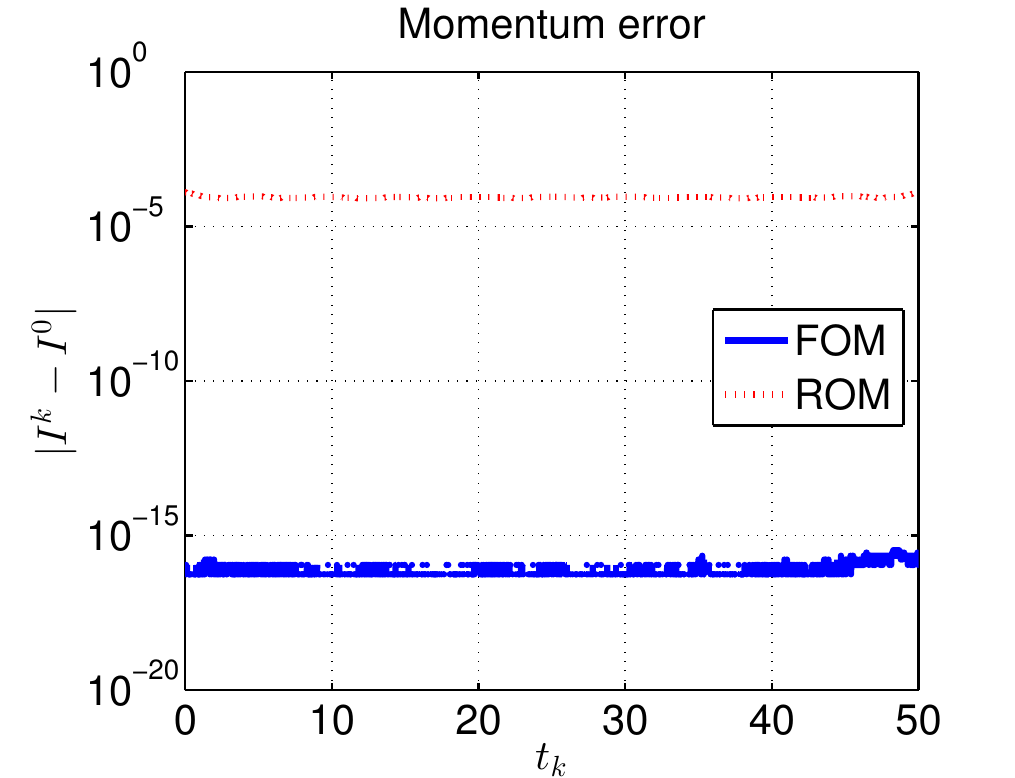}}
\caption{Preservation of the energy (left)  and momentum (right) }
\label{hamcale}
\end{figure}

\begin{figure}[H]
	\centerline{\includegraphics[width=1.4\columnwidth]{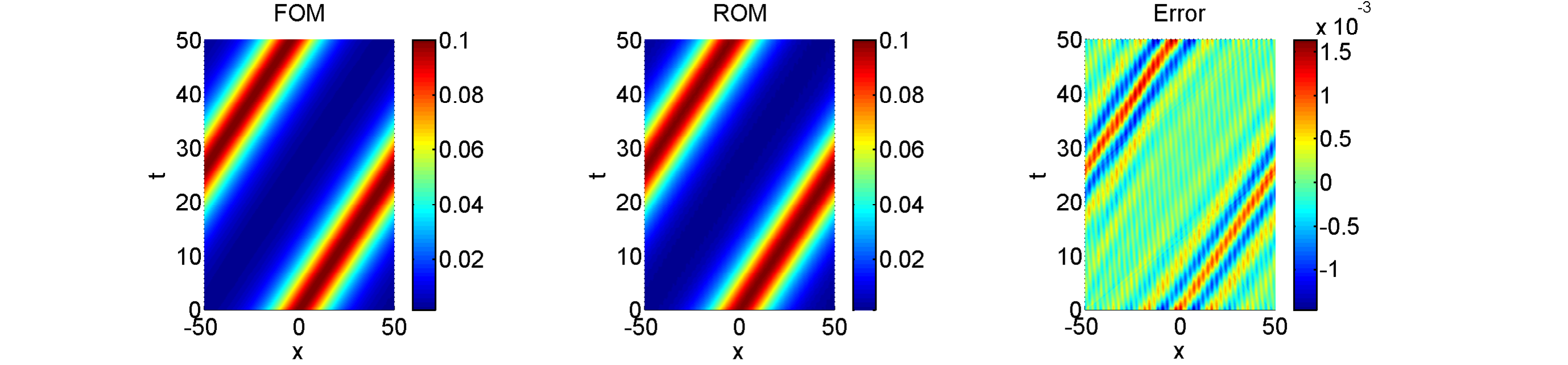}}
\caption{Solution profiles and error plot}
\label{solcale}
\end{figure}

In Figure~\ref{hamcale}, the Hamiltonian and momentum are preserved over the time without a drift, which ensures the stability of the solitons. In Table \ref{tbcale}, the
relative FOM-ROM errors \eqref{relerr} and time-averaged errors \ref{conserr1} of the invariants do not decrease much when the POD/DEIM modes are increased. Therefore, the ROM solutions in Figure \ref{solcale} accurately capture the full-order solutions.

\begin{table}[H]
\centering
\caption{Time-averaged solution and energy/momentum errors \label{tbcale}}
\begin{tabular}{|c|ccc|}\hline
$\#$ POD/DEIM Modes &   $\|\bm{\psi}-\widehat{\bm  \psi}\|_{rel}$ & $\|H\|_{\text{abs}}$ &   $\|I\|_{\text{abs}}$  \\
	\hline
    20   &   1.06e-02 &    2.95e-04  &   2.95e-04 \\ \hline
    30  &    1.56e-03  &   4.82e-06  &   3.50e-06 \\ \hline
    40  &    1.73e-03  &   3.06e-07  &   3.29e-06 \\ \hline
    50  &    1.79e-03  &   2.21e-07  &   5.70e-06  \\  \hline
\end{tabular}
\end{table}

\subsection{Damped ALE}

\begin{figure}[H]
	\centerline{\includegraphics[width=0.5\columnwidth]{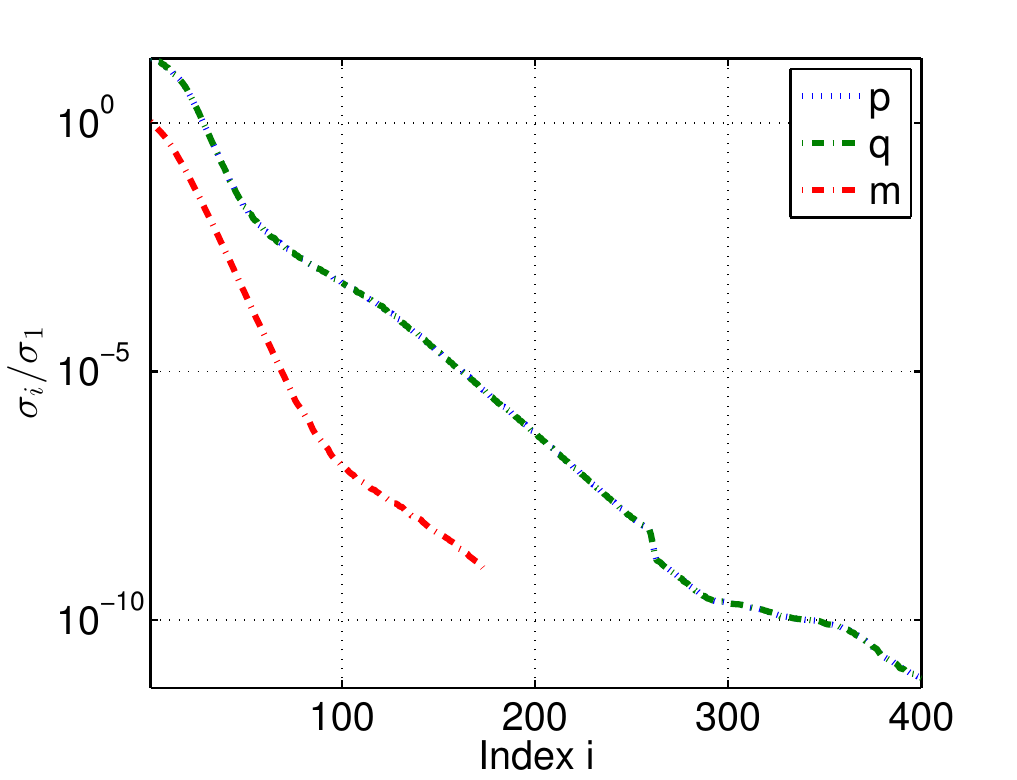}}
\caption{Singular values to the states and the nonlinearity}
\label{singdal}
\end{figure}

We consider the damped NLSE \eqref{dnls} with one soliton solution with $\gamma=0.5$ and $\mu=0.01$ \cite{Fu16},  and with the initial condition 
$$
\psi(x,0) = \frac{\sqrt{2}}{2} {\rm exp}\left ( i\frac{x+p}{2}\right ) {\rm sech} \left ( \frac{x+p}{2}\right ),
$$
where $p=20$.  We take the time step $\Delta t = 0.01$ in the time interval $t\in [0,60]$, and on the spatial domain $x\in [-64,64]$ with the mesh size $h=0.25$. The resulting FOM is of size $N=512$. The full-order solutions are saved at every five time steps, for which the snapshot matrices are of size $512\times 1200$. In Figure~\ref{singdal}, normalized singular values are presented for the damped ALE.

\begin{figure}[H]
	\centerline{\includegraphics[width=1\columnwidth]{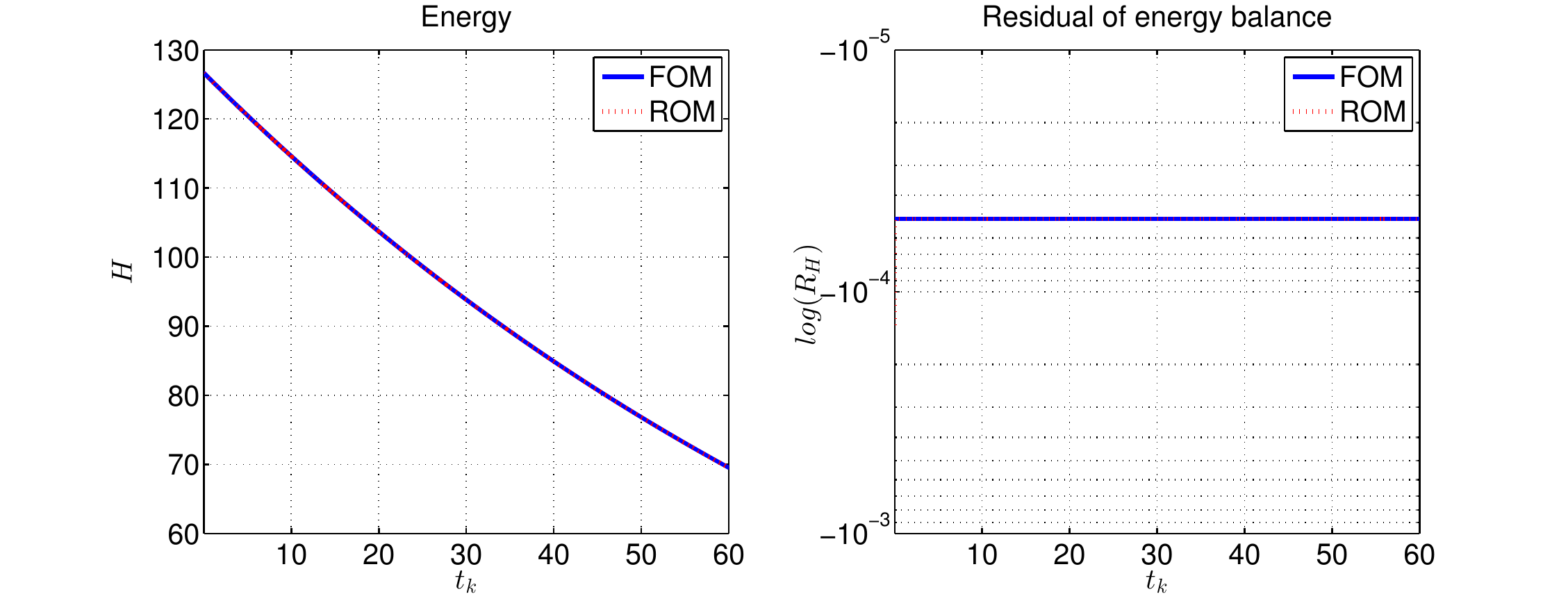}}
\caption{Energy (left) and  residual of energy balance (right) }
\label{enbalance}
\end{figure}

\begin{figure}[H]
	\centerline{\includegraphics[width=1\columnwidth]{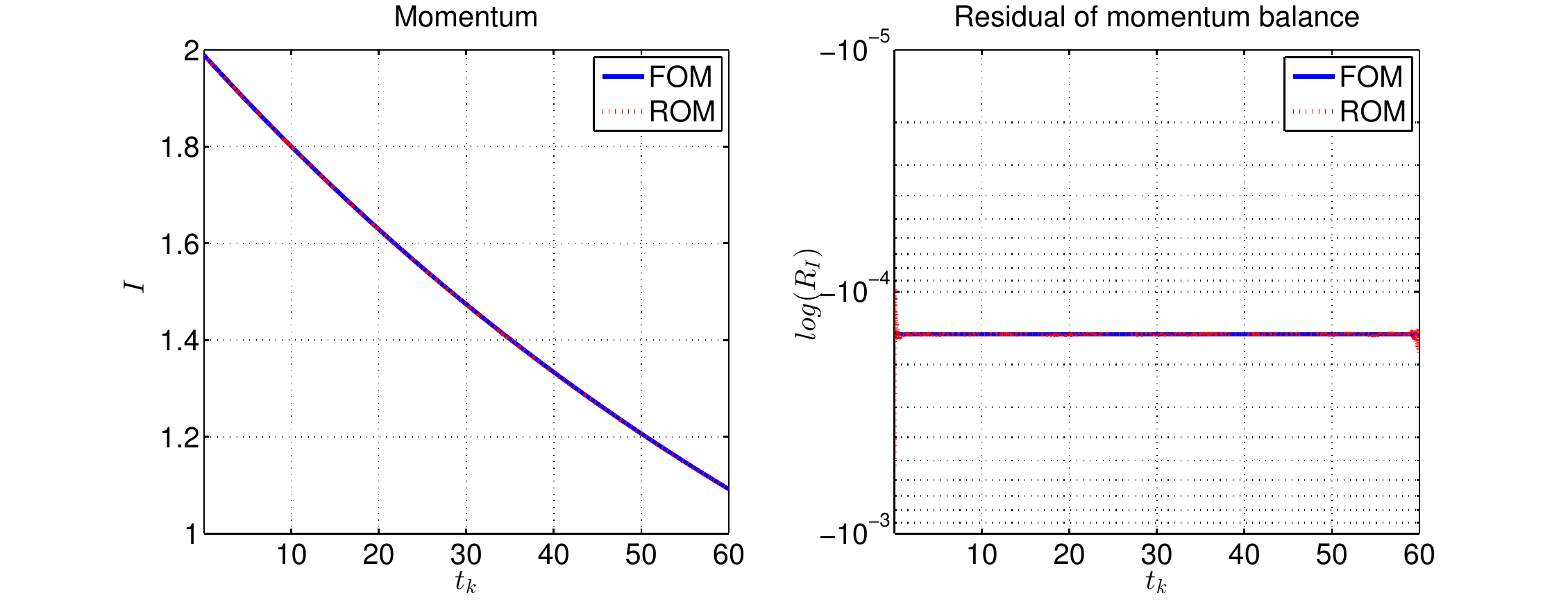}}
\caption{Momentum (left)  and  residual of momentum balance (right)}
\label{mombalance}
\end{figure}

\begin{figure}[H]
	\centerline{\includegraphics[width=0.5\columnwidth]{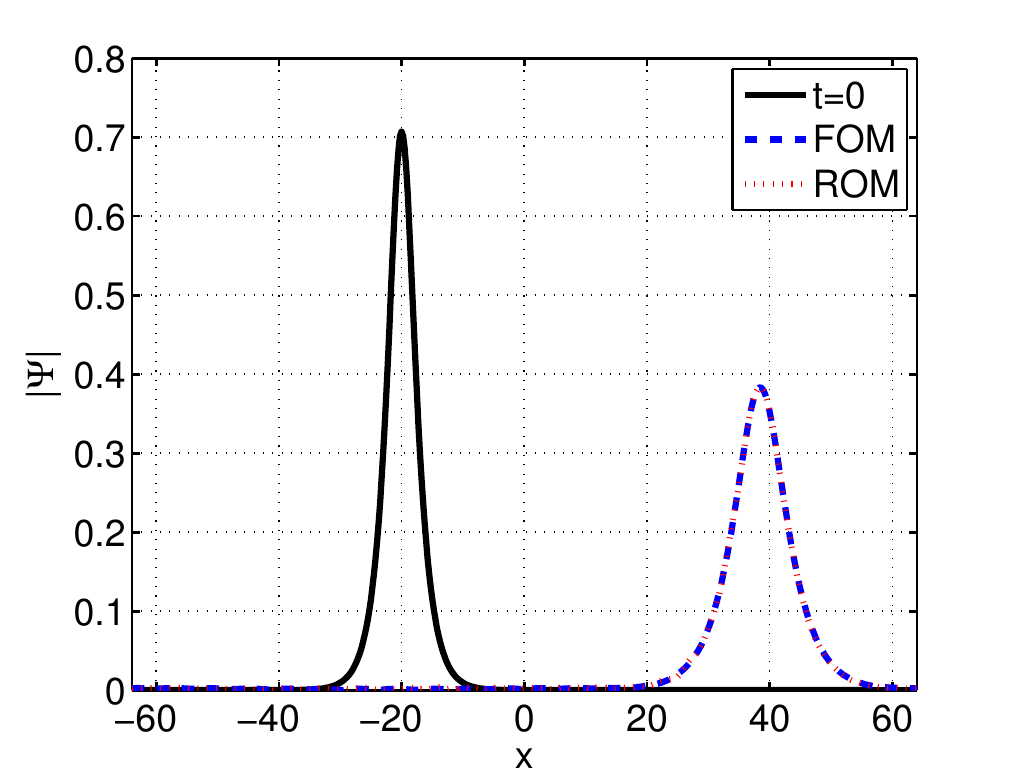}}	
\caption{Initial soliton and soliton solutions at final time}	
\label{fomromsol1_damped}
\end{figure}

In order to fulfill the energy and momentum balances in Figures \ref{enbalance}-\ref{mombalance}, the tolerances are taken as one order smaller than for the conservative ALE, because the singular values decay monotonically without reaching a plateau in contrast to the decay of singular values for the conservative ALE in Figure \ref{singal}.
By setting the tolerances $\kappa = 10^{-5}$ and $\kappa = 10^{-7}$, we obtain  $N_r=40$ POD and $N_d=49$ DEIM modes, respectively.  
The relative FOM-ROM error \eqref{relerr} is 7.09e-03, whereas the relative energy and momentum balance errors \eqref{conserr2} are 1.38-04 and 5.24e-04, respectively. The solution error and residuals of the energy/momentum balances are saturated around 50 POD/DEIM modes in Table \ref{tbdale}, similar to the conservative ALE in Figure \ref{tbcale}. The FOM/ROM solutions in Figure \ref{fomromsol1_damped} are almost identical, and the energy and momentum dissipate with correct rates in Figures \ref{enbalance}-\ref{mombalance}.

\begin{table}[H]
\centering
\caption{Time-averaged solution and energy/momentum balance errors \label{tbdale}}
\begin{tabular}{|c|ccc|}\hline
$\#$ POD/DEIM Modes &   $\|\bm{\psi}-\widehat{\bm  \psi}\|_{rel}$ & $\|R_H\|_{\text{abs}}$ &   $\|R_I\|_{\text{abs}}$  \\
	\hline
    20   &   2.24e-01  &   9.64e-02  &   1.95e-01 \\ \hline
    30   &   3.89e-02   &  3.19e-03  &   8.55e-03 \\ \hline
    40   &   7.07e-03   &  1.38e-04  &   5.25e-04 \\ \hline
    50   &   1.70e-03   &  5.39e-05  &   1.65e-04 \\ \hline
    60   &    7.23e-04  &   5.04e-05  &     1.51e-04  \\ \hline
\end{tabular}
\end{table}

\begin{figure}[H] 
	\centerline{\includegraphics[width=0.5\columnwidth]{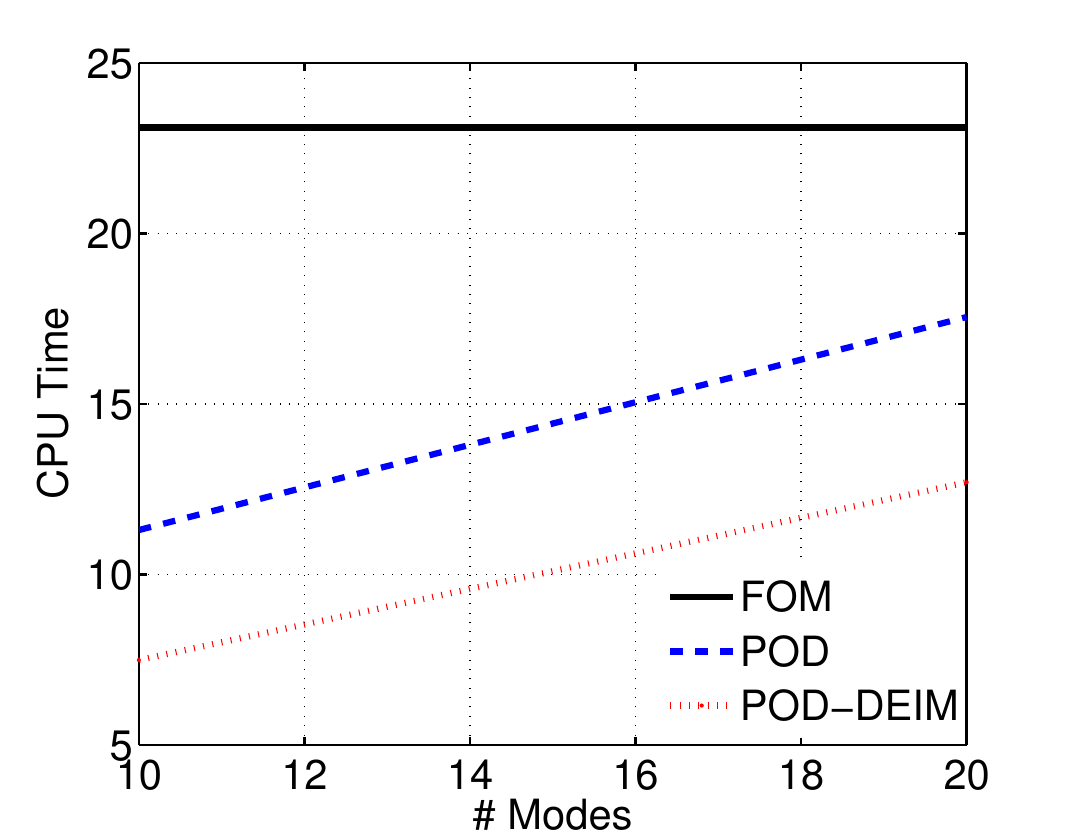}
		 \includegraphics[width=0.5\columnwidth]{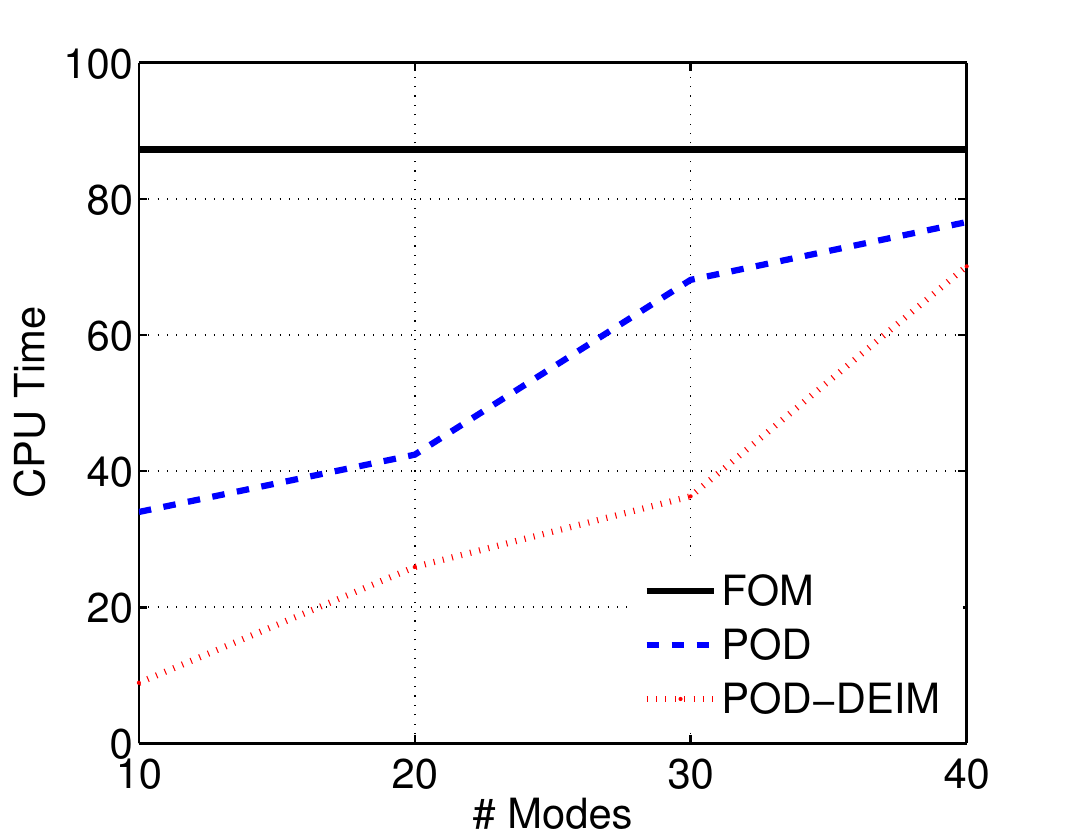}}
\caption{CPU times of FOM and different ROMs for conservative ALE (left) and damped ALE (right) }
\label{cpu}
\end{figure}

In Figure \ref{cpu}, we present the CPU times (in seconds) versus number of POD/DEIM modes for the conservative and damped ALEs. 
It is seen that the ROM becomes inefficient after $20$ modes for the conservative ALE, and $40$ modes for the damped ALE. 
For larger number of modes, ROMs may need more computing time even than the FOM. This is due to the term ${\bm m}_r\otimes {\bm p}_r$ to be computed in the POD-DEIM ROM \eqref{deim3}, whose computational cost is $\mathcal{O}(N_rN_d)$. Generally, it is expected for computational efficiency that $N_rN_d<N$, but this can be violated for large number of modes causing that $N_rN_d>N$.
The maximum number of the modes in Figure \ref{cpu} correspond to those in Figure \ref{solcale} and Figure \ref{fomromsol1_damped}, which are determined by the energy criteria \eqref{energy_criteria}. On the other side, it is not necessary to compute the ROMs beyond $20$ POD modes for the conservative ALE, because the energy and momentum errors are almost the same for all modes greater than $20$, as shown in Table \ref{tbcale}. 
Similarly, the time-averaged solution and energy/momentum balance errors do not decrease significantly after $40$ POD modes in Table \ref{tbdale} for the damped ALE. 
The effect of the tensorial ROM might be more visible for the problems on a spatial domain of two dimension or higher. 
In \cite{Karasozen21sw,Karasozen22rtswe,Stefanescu14}, the authors construct computationally efficient ROMs with tensor techniques for the two-dimensional shallow water equations that preserve the Hamiltonian, mass and other conserved quantities.

\section{Conclusions}
\label{con}

In this paper, the ROMs that preserve the conservation and dissipation structures of the conservative and dissipative ALEs, respectively, are constructed.
Preservation of structures of the conservative and dissipative ALEs in reduced-order form, ensures the stability and robustness of soliton solutions in long term. 
As common in linear MOR techniques, the POD and DEIM require relatively large number of modes to resolve the reduced-order solutions accurately, as shown in the numerical examples.
The recently developed nonlinear MOR on manifolds \cite{Buchfink21} or kernel methods \cite{Diez21} can produce more accurate solutions in low dimensional reduced spaces. As a future work, we aim to develop structure-preserving nonlinear ROMs for canonical and non-canonical Hamiltonian systems. Using the scalar auxiliary variable (SAV) approach \cite{Tapley22,Shen19,Kemmochi22}, the polynomial higher order conserved quantities of Hamiltonian PDEs, like the Korteweg de Vries equation, and energy functionals of dissipative PDEs, like the Allen-Cahn and Chan-Hilliard equations, are expressed as quadratic invariants. The transformed PDEs and ODEs can be integrated by the mid-point rule. The tensor techniques developed in this paper can be applied in a ROM framework efficiently.

\section*{Acknowledgments}
The authors thank the referees for their constructive comments which helped much to improve the paper.


\begin{thebibliography}{10}

\bibitem{Ablowitz76}
M.~J. Ablowitz and J.~F. Ladik.
\newblock A nonlinear difference scheme and inverse scattering.
\newblock {\em Studies in Applied Mathematics}, 55(3):213--229, 1976.

\bibitem{Hesthaven16}
B.~M. Afkham and J.S. Hesthaven.
\newblock Structure preserving model reduction of parametric {Hamiltonian}
  systems.
\newblock {\em SIAM Journal on Scientific Computing}, 39(6):A2616--A2644, 2017.

\bibitem{Ahmed20}
S.~E. Ahmed and O.~San.
\newblock Breaking the {Kolmogorov} barrier in model reduction of fluid flows.
\newblock {\em Fluids}, 5(1), 2020.

\bibitem{Benner21}
P.~Benner and P.~Goyal.
\newblock Interpolation-based model order reduction for polynomial systems.
\newblock {\em SIAM Journal on Scientific Computing}, 43(1):A84--A108, 2021.

\bibitem{Benner18}
P.~Benner, P.~Goyal, and S.~Gugercin.
\newblock ${\mathcal h}_2$-quasi-optimal model order reduction for
  quadratic-bilinear control systems.
\newblock {\em SIAM Journal on Matrix Analysis and Applications},
  39(2):983--1032, 2018.

\bibitem{Benner20hbmor}
P.~Benner, S.~Grivet-Talocia, A.~Quarteroni, G.~Rozza, W.~Schilders, and M.L.
  Silveira.
\newblock {\em Model Order Reduction: Volume 2: Snapshot-Based Methods and
  Algorithms}.
\newblock De Gruyter, 2020.

\bibitem{Berkooz93}
G~Berkooz, P~Holmes, and J~L Lumley.
\newblock The proper orthogonal decomposition in the analysis of turbulent
  flows.
\newblock {\em Annual Review of Fluid Mechanics}, 25(1):539--575, 1993.

\bibitem{Bhatt16}
A.~Bhatt, D.~Floyd, and B.~E. Moore.
\newblock Second order conformal symplectic schemes for damped {Hamiltonian}
  systems.
\newblock {\em Journal of Scientific Computing}, 66(3):1234--1259, 2016.

\bibitem{Buchfink18}
P.~Buchfink, A.~Bhatt, and B.~Haasdonk.
\newblock Symplectic model order reduction with non-orthonormal bases.
\newblock {\em Mathematical \& Computational Applications}, 24(2):Paper No. 43,
  26, 2019.

\bibitem{Buchfink21}
P.~Buchfink, S.~Glas, and B.~Haasdonk.
\newblock Symplectic model reduction of {Hamiltonian} systems on nonlinear
  manifolds, 2021.

\bibitem{Buchfink20greedy}
P.~Buchfink, B.~Haasdonk, and S.~Rave.
\newblock Psd-greedy basis generation for structure-preserving model order
  reduction of {Hamiltonian} systems.
\newblock {\em Proceedings of the Conference Algoritmy}, pages 151--160, 2020.

\bibitem{Carlberg15}
K.~Carlberg, R.~Tuminaro, and P.~Boggs.
\newblock Preserving {L}agrangian structure in nonlinear model reduction with
  application to structural dynamics.
\newblock {\em SIAM J. Sci. Comput.}, 37(2):B153--B184, 2015.

\bibitem{Chaturantabu16}
S.~Chaturantabut, C.~Beattie, and S.~Gugercin.
\newblock Structure-preserving model reduction for nonlinear {port-Hamiltonian}
  systems.
\newblock {\em SIAM Journal on Scientific Computing}, 38(5):B837--B865, 2016.

\bibitem{chaturantabut10nmr}
S.~Chaturantabut and D.~C. Sorensen.
\newblock Nonlinear model reduction via discrete empirical interpolation.
\newblock {\em SIAM J. Sci. Comput.}, 32(5):2737--2764, 2010.

\bibitem{Cohen11}
D.~Cohen and E.~Hairer.
\newblock Linear energy-preserving integrators for {Poisson} systems.
\newblock {\em BIT Numerical Mathematics}, 51(1):91--101, 2011.

\bibitem{Stefanescu14}
R\u{a}zvan \c{S}tef\u{a}nescu, Adrian Sandu, and Ionel~M. Navon.
\newblock Comparison of {POD} reduced order strategies for the nonlinear {2D}
  shallow water equations.
\newblock {\em International Journal for Numerical Methods in Fluids},
  76(8):497--521, 2014.

\bibitem{Diez21}
P.~D\'{\i}ez, A.~Muix\'{\i}, S.~Zlotnik, and A.~Garc\'{\i}a-Gonz\'{a}lez.
\newblock Nonlinear dimensionality reduction for parametric problems: A kernel
  proper orthogonal decomposition.
\newblock {\em International Journal for Numerical Methods in Engineering},
  122(24):7306--7327, 2021.

\bibitem{gugercin16a}
Z.~Drma\v{c} and S.~Gugercin.
\newblock A new selection operator for the discrete empirical interpolation
  method--improved a priori error bound and extensions.
\newblock {\em SIAM Journal on Scientific Computing}, 38(2):A631--A648, 2016.

\bibitem{Fu16}
H.~Fu, W.E. Zhou, X.~Qian, S.H. Song, and L.Y. Zhang.
\newblock Conformal structure-preserving method for damped nonlinear
  {Schr{\"o}dinger} equation.
\newblock {\em Chinese Physics B}, 25(11):110201, 2016.

\bibitem{Gong17}
Y.~Gong, Q.~Wang, and .~Wang.
\newblock Structure-preserving {Galerkin} {POD} reduced-order modeling of
  {H}amiltonian systems.
\newblock {\em Computer Methods in Applied Mechanics and Engineering}, 315:780
  -- 798, 2017.

\bibitem{Hesthaven22an}
J.~S. Hesthaven, C.~Pagliantini, and G.~Rozza.
\newblock Reduced basis methods for time-dependent problems.
\newblock {\em Acta Numerica}, 31:265--345, 2022.

\bibitem{Islas01}
A.L. Islas, D.A. Karpeev, and C.M. Schober.
\newblock Geometric integrators for the nonlinear {Schr\"odinger} equation.
\newblock {\em Journal of Computational Physics}, 173(1):116 -- 148, 2001.

\bibitem{Karasozen18nls}
B.~Karas\"{o}zen and M.~Uzunca.
\newblock Energy preserving model order reduction of the nonlinear
  {S}chr\"{o}dinger equation.
\newblock {\em Advances in Computational Mathematics}, 44(6):1769--1796, 2018.

\bibitem{Karasozen21sw}
B.~Karas\"{o}zen, S.~Y{\i}ld{\i}z, and M.~Uzunca.
\newblock Structure preserving model order reduction of shallow water
  equations.
\newblock {\em Mathematical Methods in the Applied Sciences}, 44(1):476--492,
  2021.

\bibitem{Karasozen22rtswe}
B.~Karas{\o}zen, S.~Y{\i}ld{\i}z, and M.~Uzunca.
\newblock Energy preserving reduced-order modeling of the rotating thermal
  shallow water equation.
\newblock {\em Physics of Fluids}, 34(5):056603, 2022.

\bibitem{Karasozen18}
B{\"u}lent Karas{\"o}zen and Murat Uzunca.
\newblock Energy preserving model order reduction of the nonlinear
  schr{\"o}dinger equation.
\newblock {\em Advances in Computational Mathematics}, 44(6):1769--1796, 2018.

\bibitem{Kemmochi22}
T.~Kemmochi and S.~Sato.
\newblock Scalar auxiliary variable approach for conservative/dissipative
  partial differential equations with unbounded energy functionals.
\newblock {\em BIT}, 62(3):903--930, 2022.

\bibitem{leva08mmm}
P.~d. Leva.
\newblock {MULTIPROD TOOLBOX}, multiple matrix multiplications, with array
  expansion enabled.
\newblock Technical report, University of Rome Foro Italico, Rome, 2008.

\bibitem{Miyatake19}
Y.~Miyatake.
\newblock Structure-preserving model reduction for dynamical systems with a
  first integral.
\newblock {\em Japan Journal of Industrial and Applied Mathematics},
  36(3):1021--1037, 2019.

\bibitem{Moore21}
B.~E. Moore.
\newblock Exponential integrators based on discrete gradients for linearly
  damped/driven {P}oisson systems.
\newblock {\em Journal of Scientific Computing}, 87(2):Paper No. 56, 18, 2021.

\bibitem{Ohlberger16}
M.~Ohlberger and S.~Rave.
\newblock Reduced basis methods: Success, limitations and future challenges.
\newblock {\em Proceedings of the Conference Algoritmy}, pages 1--12, 2016.

\bibitem{Peherstorfer22}
B.~Peherstorfer.
\newblock Breaking the {Kolmogorov} barrier with nonlinear model reduction.
\newblock {\em Notices of the American Mathematical Society}, 65(9):725--733,
  2022.

\bibitem{Peng16}
L.~Peng and K.~Mohseni.
\newblock Symplectic model reduction of {Hamiltonian} systems.
\newblock {\em SIAM Journal on Scientific Computing}, 38(1):A1--A27, 2016.

\bibitem{Schober99}
C.~M. Schober.
\newblock Symplectic integrators for the {A}blowitz-{L}adik discrete nonlinear
  {S}chr\"{o}dinger equation.
\newblock {\em Physics Letters. A}, 259(2):140--151, 1999.

\bibitem{Shen19}
J.~Shen, J.~Xu, and J.~Yang.
\newblock A new class of efficient and robust energy stable schemes for
  gradient flows.
\newblock {\em SIAM Review}, 61(3):474--506, 2019.

\bibitem{Sirovich87}
L.~Sirovich.
\newblock Turbulence and the dynamics of coherent structures. {III}. {D}ynamics
  and scaling.
\newblock {\em Quart. Appl. Math.}, 45(3):583--590, 1987.

\bibitem{Tang07}
Y.~Tang, J.~Cao, X.~Liu, and Y.~Sun.
\newblock Symplectic methods for the{ Ablowitz--Ladik }discrete nonlinear
  {Schr{\"o}dinger} equation.
\newblock {\em Journal of Physics A: Mathematical and Theoretical},
  40(10):2425--2437, 2007.

\bibitem{Tapley22}
B.~K. Tapley.
\newblock Geometric integration of {ODE}s using multiple quadratic auxiliary
  variables.
\newblock {\em SIAM Journal on Scientific Computing}, 44(4):A2651--A2668, 2022.

\bibitem{Karasozen21kdv}
M.~Uzunca, B.~Karas{\"o}zen, and S.~Y{\i}ld{\i}z.
\newblock Structure-preserving reduced-order modeling of {Korteweg-de Vries}
  equation.
\newblock {\em Mathematics and Computers in Simulation}, 188:193--211, 2021.

\bibitem{Zaharov74}
V.~E. Zaharov and S.~V. Manakov.
\newblock The complete integrability of the nonlinear {S}chr\"{o}dinger
  equation.
\newblock {\em Theoretical and Mathematical Physics}, 19:332--343, 1974.

\end{thebibliography}

\end{document}